\documentclass{amsart}

\usepackage{amssymb}
\usepackage[mathcal]{eucal}
\usepackage{upref}

\usepackage{hyperref}

\usepackage[lite]{amsrefs}

\usepackage[all,cmtip]{xy}

\newtheoremstyle{slplain}
 {.5\baselineskip\@plus.2\baselineskip\@minus.2\baselineskip}
 {.5\baselineskip\@plus.2\baselineskip\@minus.2\baselineskip}
 {\slshape}
 {}
 {\bfseries}
 {.}
 { }
 {}

\numberwithin{equation}{section}
\theoremstyle{slplain}
\newtheorem{thm}[equation]{Theorem}  
\newtheorem{cor}[equation]{Corollary}     
\newtheorem{lem}[equation]{Lemma}         
\newtheorem{prop}[equation]{Proposition}  

\theoremstyle{definition}
\newtheorem{defn}[equation]{Definition} 

\theoremstyle{remark}
\newtheorem{rem}[equation]{Remark}      
\newtheorem{ex}[equation]{Example}      

\newcommand{\thmref}{Theorem~\ref}
\newcommand{\propref}{Proposition~\ref}
\newcommand{\lemref}{Lemma~\ref}
\newcommand{\defref}{Definition~\ref}

\newcommand{\diagref}{Diagram~\ref}

\newcommand{\secref}{Section~\ref}

\newcommand{\cat}[1]{\mathcal{#1}}

\newcommand{\R}{\mathbb{R}}
\newcommand{\suchthat}{\mid}
\newcommand{\rest}[1]{|_{#1}}

\newcommand{\pr}{\mathrm{pr}}

\newcommand{\bdry}{\partial}

\DeclareMathOperator*{\colim}{colim}
\DeclareMathOperator{\Map}{Map}
\newcommand{\homotopic}{\simeq}

\newcommand{\CofSquares}{%
  \mappullback{(S^{n-1}, W_{k})}
  {(S^{n-1}, Y)}
  {(D^{n}, Y)}%
}

\newcommand{\pullback}[3]{#1\mathbin{\mathord{\times}_{\!#2}}#3}
\newcommand{\mappullback}[3]{\pullback{\Map#1}{\Map#2}{\Map#3}}

\newcommand{\Comma}{\rlap{\enspace ,}}
\newcommand{\Period}{\rlap{\enspace .}}
\newcommand{\Semicolon}{\rlap{\enspace ;}}

\begin{document}

\title[The Quillen model category of topological spaces]{The Quillen
  model category\\of topological spaces}

\author{Philip S. Hirschhorn}

\address{Department of Mathematics\\
  Wellesley College\\
  106 Central Street\\
  Wellesley, Massachusetts 02481}

\email{psh@math.mit.edu}

\urladdr{http://www-math.mit.edu/\textasciitilde psh}

\begin{abstract}
  We give a complete and careful proof of Quillen's theorem on the
  existence of the standard model category structure on the category
  of topological spaces. We do not assume any familiarity with model
  categories.
\end{abstract}

\thanks{\textcopyright 2017. This manuscript version is made available
  under the CC-BY-NC-ND 4.0 license
  \url{http://creativecommons.org/licenses/by-nc-nd/4.0/}}



\date{September 11, 2017}

\maketitle

\tableofcontents

\section{Introduction}
\label{sec:Intro}

Quillen defined model categories (see \defref{def:ModCat}) in
\cites{Q:HomAlg,Q:RHT} to apply the techniques of homotopy theory to
categories other than topological spaces or simplicial sets.  Of
course, one of his main examples was a model category structure on the
category of topological spaces, but the proof that this example
satisfied the axioms was only sketched in places, and filling in the
details is actually rather hard.  Hovey \cite{Hovey:MC} presents a
complete proof, but uses some difficult arguments involving function
spaces and adjointness.  In this note we assume no previous knowledge
of model categories and we present a complete statement and proof of
Quillen's theorem.  We assume nothing more sophisticated than homotopy
groups and the exact homotopy sequence of a fibration.

In early work in homotopy theory, the \emph{homotopy category} was
defined to be a category in which the morphisms were homotopy classes
of maps, and so the isomorphisms were the morphisms represented by
homotopy equivalences.  Weak homotopy equivalences (i.e., maps that
induced an isomorphism of the set of path components and of all
homotopy groups at all basepoints) induced isomorphisms of all
homology and cohomology groups, but they were not in general homotopy
equivalences unless both the domain and target were CW-complexes.  For
this reason, the \emph{homotopy category} was often defined to have as
objects only the CW-complexes, or spaces homotopy equivalent to a
CW-complex.  Whenever a construction led to a space that was not
homotopy equivalent to a CW-complex, it was replaced by a weakly
equivalent CW-complex; this preserved all the relevant algebraic
invariants.  Thus, weakly equivalent spaces were recognized as being
somehow ``equivalent'', even if that was never made explicit.  The
\emph{homotopy category of a model category} (see
\cite{MCATL}*{Def.~8.3.2}) is constructed by making the weak
equivalences into isomorphisms.  This identifies homotopic maps (see
\cite{MCATL}*{Lemma~8.3.4}) and takes homotopy equivalences into
isomorphisms, but it also explicitly takes all weak equivalences into
isomorphisms.

Early work considered the fibrations that are now called
\emph{Hurewicz fibrations}, maps with the homotopy lifting property
with respect to all spaces.  When Serre \cite{Serre} was developing
his spectral sequence for the homology of fiber spaces, he realized
that a weaker condition, that of having the homotopy lifting property
with respect to polyhedra (or, equivalently, with respect to cubes),
was sufficient for work involving homotopy and homology groups, and
this larger class of \emph{Serre fibrations} is now more commonly
used.  Quillen's model category of topological spaces takes as
\emph{weak equivalences} the weak homotopy equivalences, as
\emph{fibrations} the Serre fibrations, and as \emph{cofibrations} the
relative cell complexes (see \defref{def:CellComp}) and their
retracts.

There is also another model category structure on the category of
topological spaces.  Str{\o}m \cite{Strom:HC} has shown that there is
a model category structure on the category of topological spaces in
which the weak equivalences are the homotopy equivalences, the
fibrations are the Hurewicz fibrations, and the cofibrations are the
closed inclusions with the homotopy extension property.

\section{Definitions and the main theorem}
\label{sec:MCdef}

This note can be read as applying to any of the standard complete and
cocomplete categories (i.e., a category that contains both a limit and
a colimit for every small diagram in the category) of topological
spaces, e.g.,
\begin{itemize}
\item the category of all topological spaces,
\item the category of compactly generated spaces (see, e.g.,
  \cite{Boardman-Vogt}*{Appendix}),
\item the category of compactly generated weak Hausdorff spaces (see,
  e.g., \cite{Strickland}, or \cite{Fritsch-Piccinini}*{Appendix}), or
\item the category of compactly generated Hausdorff spaces (see, e.g.,
  \cite{Steenrod}).
\end{itemize}
It also applies to any other category of spaces that contains all
inclusions of a subcomplex of a CW-complex and in which
\propref{prop:RelCellCompSubset} (that a compact subset of a relative
cell complex intersects the interiors of only finitely many cells)
holds.  The main theorem (on the existence of the model category
structure) is \thmref{thm:main}.

\begin{defn}
  \label{def:retract}
  If there is a commutative diagram
  \begin{displaymath}
    \xymatrix{
      {A} \ar[r] \ar[d]_{i} 
           \ar `u[rr] `[drr]^{1_{A}} [rr]+<0pt,6pt>
        & {C} \ar[r] \ar[d]^{j}
        & {A} \ar[d]^{i}\\
      {B} \ar[r]
           \ar `d[rr] `[urr]_{1_{B}} [rr]
        & {D} \ar[r]
        & {B}
    }
  \end{displaymath}
  then we will say that the map $i$ is a \emph{retract} of the map
  $j$.
\end{defn}

We use the definition of model category from
\cite{MCATL}*{Def.~7.1.3}; this differs from Quillen's original
definition of a \emph{closed model category} (\cites{Q:HomAlg,Q:RHT})
in that we require the existence of \emph{all} colimits and limits
(not just the finite ones), and we require that the two factorizations
be functorial.

\begin{defn}
  \label{def:ModCat}
  A \emph{model category} is a category $\cat M$ together with three
  classes of maps (called the \emph{weak equivalences}, the
  \emph{cofibrations}, and the \emph{fibrations}), satisfying the
  following five axioms:
  \begin{itemize}
  \item[M1:] (Limit axiom) The category $\cat M$ is complete and
    cocomplete (i.e., contains both a limit and a colimit for every
    small diagram in the category).
  \item[M2:] (Two out of three axiom) If $f$ and $g$ are maps in $\cat
    M$ such that $gf$ is defined and two of $f$, $g$, and $gf$ are
    weak equivalences, then so is the third.
  \item[M3:] (Retract axiom) If $f$ and $g$ are maps in $\cat M$ such
    that $f$ is a retract of $g$ (in the category of maps of $\cat M$;
    see \defref{def:retract}) and $g$ is a weak equivalence, a
    fibration, or a cofibration, then so is $f$.
  \item[M4:] (Lifting axiom) Given the commutative solid arrow diagram
    in $\cat M$
    \begin{displaymath}
      \xymatrix{
        {A} \ar[r] \ar[d]_{i} 
          & {X} \ar[d]^{p}\\
        {B} \ar@{.>}[ur] \ar[r]  
          & {Y}
      }
    \end{displaymath}
    the dotted arrow exists if either
    \begin{enumerate}
    \item $i$ is a cofibration and $p$ is both a fibration and a weak
      equivalence or
    \item $i$ is both a cofibration and a weak equivalence and $p$ is
      a fibration.
    \end{enumerate}
  \item[M5:] (Factorization axiom) Every map $g$ in $\cat M$ has two
    functorial factorizations:
    \begin{enumerate}
    \item $g = pi$, where $i$ is a cofibration and $p$ is both a
      fibration and a weak equivalence, and
    \item $g = qj$, where $j$ is both a cofibration and a weak
      equivalence and $q$ is a fibration.
    \end{enumerate}
  \end{itemize}
\end{defn}

\begin{defn}
  \label{def:TrivCof}
  Let $\cat M$ be a model category.
  \begin{enumerate}
  \item A \emph{trivial fibration} is a map that is both a fibration
    and a weak equivalence.
  \item A \emph{trivial cofibration} is a map that is both a
    cofibration and a weak equivalence.
  \item An object is \emph{cofibrant} if the map to it from the
    initial object is a cofibration.
  \item An object is \emph{fibrant} if the map from it to the terminal
    object is a fibration.
  \item An object is \emph{cofibrant-fibrant} if it is both cofibrant
    and fibrant.
  \end{enumerate}
\end{defn}

\begin{defn}
  \label{def:MdCtStr}
  We will say that a map $f\colon X \to Y$ of topological spaces is
  \begin{itemize}
  \item a \emph{weak equivalence} if it is a weak homotopy
    equivalence, i.e., if either $X$ and $Y$ are both empty, or $X$
    and $Y$ are both nonempty and for every choice of basepoint $x \in
    X$ the induced map $f_{*}\colon \pi_{i}(X,x) \to \pi_{i}(Y, f(x))$
    of homotopy groups (if $i > 0$) or sets (if $i = 0$) is an
    isomorphism,
  \item a \emph{cofibration} if it is a relative cell complex (see
    \defref{def:CellComp}) or a retract (see \defref{def:retract}) of
    a relative cell complex, and
  \item a \emph{fibration} if it is a Serre fibration.
  \end{itemize}
\end{defn}

\begin{thm}
  \label{thm:main}
  There is a model category structure on the category of topological
  spaces in which the weak equivalences, cofibrations, and fibrations
  are as in \defref{def:MdCtStr}.
\end{thm}

The proof of \thmref{thm:main} is in \secref{sec:proofmain}.

\section{Lifting}
\label{sec:retracts}

\begin{defn}
  \label{def:LLP}
  If $i\colon A \to B$ and $f\colon X \to Y$ are maps such that for
  every solid arrow diagram
  \begin{displaymath}
    \xymatrix{
      {A} \ar[r] \ar[d]_{i}
      & {X} \ar[d]^{f}\\
      {B} \ar[r] \ar@{..>}[ur]
      & {Y}
    }
  \end{displaymath}
  there exists a diagonal arrow making both triangles commute, then
  $i$ is said to have the \emph{left lifting property} with respect to
  $f$ and $f$ is said to have the \emph{right lifting property} with
  respect to $i$.
\end{defn}

Thus, the lifting axiom M4 (see \defref{def:ModCat}) asserts that
\begin{itemize}
\item cofibrations have the left lifting property with respect to all
  trivial fibrations, and
\item fibrations have the right lifting property with respect to all
  trivial cofibrations.
\end{itemize}
In this section, we show that the class of maps with the left lifting
property with respect to a map $f\colon X \to Y$ is closed under
pushouts (see \lemref{lem:pushLLP}), coproducts (see
\lemref{lem:coprodLLP}), transfinite compositions (see
\lemref{lem:compLLP}) and retracts (see \lemref{lem:retractLLP}).
These results are combined in \propref{prop:ClCmpLift} to show that if
a map has the right lifting property with respect to the inclusions
$S^{n-1} \to D^{n}$ for all $n \ge 0$ then it has the right lifting
property with respect to all relative cell complexes (see
\defref{def:CellComp}) and their retracts.  This is a key step in
proving that the lifting axiom M4 is satisfied (see the proof of
\thmref{thm:onelift}).  We also present the \emph{retract argument}
(see \propref{prop:RetArg}), which is often used to show that a map
has the (left or right) lifting property with respect to other maps
(see the proof of \thmref{thm:otherlift}).

\subsection{Pushouts, pullbacks, and lifting}
\label{sec:PushPullLift}

\begin{defn}
  \label{def:pushout}
  \leavevmode
  \begin{enumerate}
  \item If there is a pushout diagram
    \begin{displaymath}
      \xymatrix{
        {A} \ar[r] \ar[d]_{i}
        & {C} \ar[d]^{j}\\
        {B} \ar[r]
        & {D}
      }
    \end{displaymath}
    then we will say that the map $j$ is a \emph{pushout} of the map
    $i$.
  \item If there is a pullback diagram
    \begin{displaymath}
      \xymatrix{
        {W} \ar[r] \ar[d]_{g}
        & {Y} \ar[d]^{f}\\
        {X} \ar[r]
        & {Z}
      }
    \end{displaymath}
    then we will say that the map $g$ is a \emph{pullback} of the map
    $f$.
  \end{enumerate}
\end{defn}

\begin{lem}
  \label{lem:pushLLP}
  \leavevmode
  \begin{enumerate}
  \item If the map $j$ is a pushout of the map $i$ (see
    \defref{def:pushout}) and if $i$ has the left lifting property
    with respect to a map $f\colon X \to Y$, then $j$ has the left
    lifting property with respect to $f$.
  \item If the map $g$ is a pullback of the map $f$ and if $f$ has the
    right lifting property with respect to a map $i$, then $g$ has the
    right lifting property with respect to $i$.
  \end{enumerate}
\end{lem}

\begin{proof}
  We will prove part~1; the proof of part~2 is similar.

  Suppose that we have the solid arrow diagram
  \begin{displaymath}
    \xymatrix{
      {A} \ar[r]^{s} \ar[d]_{i}
      & {C} \ar[r]^{t} \ar[d]_(.3){j}
      & {X} \ar[d]^{f}\\
      {B} \ar[r]_{u} \ar@{..>}[urr]^(.3){w}
      & {D} \ar[r]_{v} \ar@{..>}[ur]_{g}
      & {Y}
    }
  \end{displaymath}
  in which the left hand square is a pushout, so that $j$ is a pushout
  of $i$.  Since $i$ has the left lifting property with respect to
  $f$, there is a map $w\colon B \to X$ such that $wi = ts$ and $fw =
  vu$.  Since the left hand square is a pushout, this induces a map
  $g\colon D \to X$ such that $gu = w$ and $gj = t$.  Since $(fg)u =
  fw = (v)u$ and $(fg)j = ft = (v)j$, the universal property of the
  pushout now implies that $fg = v$, and so $j$ has the left lifting
  property with respect to $f$.
\end{proof}

\subsection{Coproducts, transfinite composition, and lifting}
\label{sec:CoprodLift}

\begin{lem}
  \label{lem:coprodLLP}
  Let $f\colon X \to Y$ be a map.  If $S$ is a set and for every $s
  \in S$ we have a map $A_{s} \to B_{s}$ that has the left lifting
  property with respect to $f$, then the coproduct $\coprod_{s\in S}
  A_{s} \to \coprod_{s\in S} B_{s}$ has the left lifting property with
  respect to $f$.
\end{lem}

\begin{proof}
  Given the solid arrow diagram
  \begin{displaymath}
    \xymatrix{
      {\coprod_{s\in S} A_{s}} \ar[r] \ar[d]
      & {X} \ar[d]^{f}\\
      {\coprod_{s\in S} B_{s}} \ar[r] \ar@{..>}[ur]
      & {Y}
    }
  \end{displaymath}
  the diagonal arrow can be chosen on each summand $B_{s}$, and these
  together define it on the coproduct.
\end{proof}

\begin{lem}
  \label{lem:compLLP}
  Let $f\colon X \to Y$ be a map.  If
  \begin{displaymath}
    A_{0} \to A_{1} \to A_{2} \to \cdots
  \end{displaymath}
  is a sequence of maps such that $A_{n} \to A_{n+1}$ has the left
  lifting property with respect to $f\colon X \to Y$ for all $n \ge
  0$, then the map $A_{0} \to \colim_{n\ge 0} A_{n}$ has the left
  lifting property with respect to $f$.
\end{lem}

\begin{proof}
  Given the solid arrow diagram
  \begin{displaymath}
    \xymatrix@R=3ex{
      {A_{0}} \ar[r] \ar[d]
      & {X} \ar[ddd]^{f}\\
      {A_{1}} \ar[d] \ar@{..>}[ur]\\
      {\vdots} \ar[d]\\
      {\colim_{n\ge 0} A_{n}} \ar[r] \ar@{..>}[uuur]
      & {Y}
    }
  \end{displaymath}
  the diagonal arrow can be chosen inductively on each $A_{n}$ and
  these combine to define it on the colimit.
\end{proof}

\subsection{Retracts and lifting}
\label{sec:RetLft}

\begin{lem}[Retracts and lifting]
  \label{lem:retractLLP}
  \leavevmode
  \begin{enumerate}
  \item If the map $i$ is a retract of the map $j$ (see
    \defref{def:retract}) and $j$ has the left lifting property
    with respect to a map $f\colon X \to Y$, then $i$ has the left
    lifting property with respect to $f$.
  \item If the map $f$ is a retract of the map $g$ and $g$ has the
    right lifting property with respect to a map $i\colon A \to B$,
    then $f$ has the right lifting property with respect to $i$.
  \end{enumerate}
\end{lem}

\begin{proof}
  We will prove part~1; the proof of part~2 is similar.

  Suppose that we have the solid arrow diagram
  \begin{displaymath}
    \xymatrix{
      {A} \ar[r]^{p} \ar[d]_{i} 
           \ar `u[rr] `[drr]^{1_{A}} [rr]+<0pt,6pt>
      & {C} \ar[r]^{r} \ar[d]_{j}
      & {A} \ar[r]^{t} \ar[d]_(.3){i}
      & {X} \ar[d]^{f}\\
      {B} \ar[r]_{q} 
           \ar `d[rr] `[urr]_{1_{B}} [rr]
      & {D} \ar[r]_{s} \ar@{..>}[urr]^(.3){v}
      & {B} \ar[r]_{u} \ar@{..>}[ur]_{w}
      & {Y \Period}
    }
  \end{displaymath}
  Since $j$ has the left lifting property with respect to $f$, there
  exists a map $v\colon D \to X$ such that $vj = tr$ and $fv = us$.
  We define $w\colon B \to X$ as $w = vq$.  We then have $wi = vqi =
  vjp = trp = t$ and $fw = fvq = usq = u$.
\end{proof}

The next result, known as the \emph{retract argument}, is often used
to show that a map has the (left or right) lifting property with
respect to other maps.  This is because \lemref{lem:retractLLP} shows
that if one map is a retract of another, and if the second map has the
(left or right) lifting property with respect to another map, so does
the first.  We will use \propref{prop:RetArg} in \secref{sec:lift} to
prove that the lifting axiom M4 is satisfied (see the proof of
\thmref{thm:otherlift}).
\begin{prop}[The retract argument]
  \label{prop:RetArg}
\leavevmode
  \begin{enumerate}
  \item If the map $g$ can be factored as $g = pi$ where $g$ has the
    left lifting property with respect to $p$, then $g$ is a retract
    of $i$.
  \item If the map $g$ can be factored as $g = pi$ where $g$ has the
    right lifting property with respect to $i$, then $g$ is a retract
    of $p$.
  \end{enumerate}
\end{prop}

\begin{proof}
  We will prove part~1; the proof of part~2 is dual.

  We have the solid arrow diagram
  \begin{displaymath}
    \xymatrix{
      {X} \ar[r]^{i} \ar[d]_{g}
        & {Z} \ar[d]^{p}\\
      {Y} \ar@{.>}[ur]^{q} \ar@{=}[r]
        & {Y \Period}
    }
  \end{displaymath}
  Since $g$ has the left lifting property with respect to $p$, the
  dotted arrow $q$ exists.  This yields the commutative diagram
  \begin{displaymath}
    \xymatrix{
      {X} \ar@{=}[r] \ar[d]_{g}
        & {X} \ar@{=}[r] \ar[d]_{i}
        & {X} \ar[d]^{g}\\
        {Y} \ar[r]_{q} \ar `d[rr] `[urr]_{1_{Y}} [rr]+<0ex,-1ex>
        & {Z} \ar[r]_{p}
        & {Y \smash{\Comma}}
    }
  \end{displaymath}
  and so $g$ is a retract of $i$.
\end{proof}

\section{Relative cell complexes}
\label{sec:RelClCmp}

In this section we define and study \emph{relative cell complexes}
(see \defref{def:CellComp}) which, together with their retracts (see
\defref{def:retract}), are the cofibrations of the model category.

\begin{defn}
  \label{def:AttCell}
  If $X$ is a subspace of $Y$ such that there is a pushout square
  \begin{displaymath}
    \xymatrix{
      {S^{n-1}} \ar[r] \ar[d]
      & {X} \ar[d]\\
      {D^{n}} \ar[r]
      & {Y}
    }
  \end{displaymath}
  for some $n \ge 0$, then we will say that \emph{$Y$ is obtained from
  $X$ by attaching a cell}.
\end{defn}

\begin{defn}
  \label{def:CellComp}
  A \emph{relative cell complex} is an inclusion of a subspace
  $f\colon X \to Y$ such that $Y$ can be constructed from $X$ by a
  (possibly infinite) process of repeatedly attaching cells (see
  \defref{def:AttCell}), and it is a \emph{finite relative cell
    complex} if it can be constructed by attaching finitely many
  cells.  The topological space $X$ is a \emph{cell complex} if the
  map $\emptyset \to X$ is a relative cell complex, and it is a
  \emph{finite cell complex} if $X$ can be constructed from
  $\emptyset$ by attaching finitely many cells.
\end{defn}

\begin{ex}
  \label{ex:CWcomp}
  Every relative CW-complex is a relative cell complex, and every
  CW-complex is a cell complex.  Since the attaching map of a cell in
  a cell complex is not required to factor through the union of lower
  dimensional cells, not all cell complexes are CW-complexes.
\end{ex}

\begin{rem}
  \label{rem:MultpleCells}
  We will often construct a relative cell complex by attaching more
  than one cell at a time.  That is, given a space $X_{0}$, a set $S$,
  and for each $s \in S$ a map $S^{(n_{s}-1)} \to X_{0}$, we may
  construct a pushout
  \begin{displaymath}
    \xymatrix{
      {\coprod_{s\in S} S^{(n_{s}-1})} \ar[r] \ar[d]
      & {X_{0}} \ar[d]\\
      {\coprod_{s\in S} D^{n_{s}}} \ar[r]
      & {X_{1}}
    }
  \end{displaymath}
  and then perform a similar construction with $X_{1}$, repeating a
  possibly infinite number of times.
\end{rem}

\begin{rem}
  \label{rem:CWcell}
  While a CW-complex can be built by a countable process of attaching
  coproducts of cells, a general cell complex may require an
  arbitrarily long transfinite construction.  This is because the
  attaching map of a cell in a cell complex is not required to factor
  through the union of lower dimensional cells.
\end{rem}

\begin{rem}
  \label{rem:CelCompPres}
  \defref{def:CellComp} implies that a relative cell complex is a map
  that can be constructed as a transfinite composition of pushouts of
  inclusions of the boundary of a cell into that cell, but there will
  generally be many different possible such constructions.  When
  dealing with a topological space that is a cell complex or a map
  that is a relative cell complex, we will often assume that we have
  chosen some specific such construction (see
  \defref{def:presentation}).  Furthermore, we may choose a
  construction of the map as a transfinite composition of pushouts of
  \emph{coproducts} of cells, i.e., we will consider constructions as
  transfinite compositions in which more than one cell is attached at
  a time.
\end{rem}

\begin{defn}
  \label{def:ordinal}
  We adopt the definition of the ordinals in which an ordinal is the
  well ordered set of all lesser ordinals, and every well ordered set
  is isomorphic to a unique ordinal (see \cite{MCATL}*{Sec.~10.1.1}).
  We will often view an ordinal as a small category with objects equal
  to the elements of the ordinal and a single map from $\alpha$ to
  $\beta$ when $\alpha \le \beta$.
\end{defn}

\begin{defn}
  \label{def:lambdaseq}
  If $C$ is a class of maps and $\lambda$ is an ordinal, then a
  \emph{$\lambda$-sequence in $C$} is a functor $X\colon \lambda \to
  \mathrm{Top}$
  \begin{displaymath}
    X_{0} \to X_{1} \to X_{2} \to \cdots X_{\beta} \to \cdots
    \qquad (\beta < \lambda)
  \end{displaymath}
  such that the map $X_{\beta} \to X_{\beta+1}$ is in $C$ for $\beta+1
  < \lambda$ and for every limit ordinal $\gamma < \lambda$ the
  induced map $\colim_{\beta<\gamma} X_{\beta} \to X_{\gamma}$ is an
  isomorphism.  Equivalently, a $\lambda$-sequence in $\cat C$ is a
  colimit-preserving functor from $\lambda$ to $\mathrm{Top}$ that
  takes every map $\beta \to \beta+1$ to an element of $\cat C$.  The
  \emph{composition} of the $\lambda$-sequence is the map $X_{0} \to
  \colim_{\beta<\lambda}X_{\beta}$ (see \cite{MCATL}*{Sec.~10.2}).
\end{defn}

\begin{defn}
  \label{def:presentation}
  If $f\colon X \to Y$ is a relative cell complex, then a
  \emph{presentation} of $f$ is a $\lambda$-sequence of pushouts of
  coproducts of elements of $I$
  \begin{displaymath}
    X_{0} \to X_{1} \to X_{2} \to \cdots X_{\beta} \to \cdots
    \qquad (\beta < \lambda)
  \end{displaymath}
  such that the composition $X_{0} \to \colim_{\beta <
    \lambda}X_{\beta}$ is isomorphic to $f$.  If $e$ is a cell of the
  relative cell complex, then the \emph{presentation ordinal} of $e$
  is the first ordinal $\beta$ such that $e$ is in $X_{\beta}$.
\end{defn}

\subsection{Compact subsets of relative cell complexes}
\label{sec:CmpctSubset}

The main result of this section is
Corollary~\ref{cor:CellGenFiniteSubcmp}, which will be used in the
construction of the factorizations required by axiom M5 (see the proof
of \propref{prop:CofFact}).

\begin{prop}
  \label{prop:RelCellCompSubset}
  If $X \to Y$ is a relative cell complex, then a compact subset of
  $Y$ can intersect the interiors of only finitely many cells of
  $Y-X$.
\end{prop}

\begin{proof}
  Let $C$ be a compact subset of $Y$.  We construct a subset $P$ of
  $C$ by choosing one point of $C$ from the interior of each cell
  whose interior intersects $C$.  We will show that this subset $P$ of
  $C$ has no accumulation point in $C$, which implies that $P$ is
  finite, which implies that $C$ intersects the interiors of only
  finitely many cells of $Y-X$.

  Let $c \in C$; we will show that there is an open subset $U$ of $Y$
  such that $c\in U$ and $U\cap P$ is either empty or contains the one
  point $c$, which will imply that $c$ is not an accumulation point of
  $P$.
  
  Let $e_{c}$ be the unique cell of $Y-X$ that contains $c$ in its
  interior.  Since there is at most one point of $P$ in the interior
  of any cell of $Y-X$, we can choose an open subset $U_{c}$ of the
  interior of $e_{c}$ that contains no points of $P$ (except for $c$,
  if $c \in P$). We will use Zorn's lemma to show that we can enlarge
  $U_{c}$ to an open subset of $Y$ that contains no points of $P$
  (except for $c$, if $c \in P$).
  
  Let $\alpha$ be the presentation ordinal of the cell $e_{c}$.  If
  the presentation ordinal of the relative cell complex $X \to Y$ is
  $\gamma$, consider the set $T$ of ordered pairs $(\beta, U)$ where
  $\alpha \le \beta \le \gamma$ and $U$ is an open subset of $Y^\beta$
  such that $U\cap Y^\alpha = U_{c}$ and $U$ contains no points of $P$
  except possibly $c$.  We define a preorder on $T$ by defining
  $(\beta_1, U_1) < (\beta_2, U_2)$ if $\beta_1 < \beta_2$ and
  $U_2\cap Y^{\beta_1} = U_1$.

  If $\{(\beta_s, U_s)\}_{s\in S}$ is a chain in $T$, then
  $(\bigcup_{s\in S} \beta_s, \bigcup_{s\in S}U_s)$ is an upper bound in
  $T$ for the chain, and so Zorn's lemma implies that $T$ has a
  maximal element $(\beta_m, U_m)$.  We will complete the proof by
  showing that $\beta_m = \gamma$.
  
  If $\beta_m < \gamma$, then consider the cells of presentation
  ordinal $\beta_m + 1$.  Since $Y$ has the weak topology determined
  by $X$ and the cells of $Y-X$, we need only enlarge $U_m$ so that
  its intersection with each cell of presentation ordinal $\beta_m +
  1$ is open in that cell, and so that it still contains no points of
  $P$ except possibly $c$.  If $h\colon S^{n-1} \to Y^{\beta_m}$ is
  the attaching map for a cell of presentation ordinal $\beta_m + 1$,
  then $h^{-1}U_m$ is open in $S^{n-1}$, and so we can ``thicken''
  $h^{-1}U_m$ to an open subset of $D^n$, avoiding the (at most one)
  point of $P$ that is in the interior of the cell.  If we let $U'$
  equal the union of $U_m$ with these thickenings in the interiors of
  the cells of presentation ordinal $\beta_m + 1$, then the pair
  $(\beta_m + 1, U')$ is an element of $T$ greater than the maximal
  element $(\beta_m,U_m)$ of $T$.  This contradiction implies that
  $\beta_m = \gamma$, and so the proof is complete.
\end{proof}

\begin{prop}
  \label{prop:SubCount}
  Every cell of a relative cell complex is contained in a finite
  subcomplex of the relative cell complex.
\end{prop}

\begin{proof}
  Choose a presentation of the relative cell complex $X \to Y$ (see
  \defref{def:presentation}); we will prove the proposition by a
  transfinite induction on the presentation ordinal of the cell.  The
  induction is begun because the very first cells attached are each in
  a subcomplex with only one cell.  Since the presentation ordinal of
  every cell is a successor ordinal, it is sufficient to assume that
  the result is true for all cells of presentation ordinal at most
  some ordinal $\beta$ and show that it is also true for cells of
  presentation ordinal the successor of $\beta$.

  The image of the attaching map of any cell of presentation ordinal
  the successor of $\beta$ is compact, and so
  \propref{prop:RelCellCompSubset} implies that it intersects the
  interiors of only finitely many cells, each of which (by the
  induction hypothesis) is contained in a finite subcomplex.  The
  union of those finite subcomplexes and this new cell is then a
  finite subcomplex containing the new cell.
\end{proof}

\begin{cor}
  \label{cor:CellGenFiniteSubcmp}
  A compact subset of a relative cell complex is contained in a finite
  subcomplex of the relative cell complex.
\end{cor}

\begin{proof}
  \propref{prop:RelCellCompSubset} implies that a compact subset
  intersects the interiors of only finitely many cells, and
  \propref{prop:SubCount} implies that each of those cells is
  contained in a finite subcomplex; the union of those finite
  subcomplexes thus contains our compact subset.
\end{proof}

\subsection{Relative cell complexes and lifting}
\label{sec:ClCmpLift}

\begin{defn}
  \label{def:gc}
  We will let $I$ denote the set of maps
  \begin{displaymath}
    I = \bigl\{S^{n-1} \to D^{n} \suchthat n \ge 0\bigr\}
  \end{displaymath}
  and we will call $I$ the \emph{set of generating cofibrations}.
\end{defn}

\begin{prop}
  \label{prop:ClCmpLift}
  If a map has the right lifting property (see \defref{def:LLP}) with
  respect to every element of $I$ (see \defref{def:gc}), then it has
  the right lifting property with respect to all relative cell
  complexes and their retracts (see \defref{def:retract}).
\end{prop}

\begin{proof}
  If a map has the right lifting property with respect to every
  element of $I$, then \lemref{lem:coprodLLP} implies that it has the
  right lifting property with respect to every coproduct of elements
  of $I$, and so \lemref{lem:pushLLP} implies that it has the right
  lifting property with respect to every pushout of a coproduct of
  elements of $I$, and so \lemref{lem:compLLP} implies that it has the
  right lifting property with respect to every composition of pushouts
  of coproducts of elements of $I$, and so \lemref{lem:retractLLP}
  implies that it has the right lifting property with respect to every
  retract of such a composition.
\end{proof}

\section{The small object argument}
\label{sec:SmObj}

In this section we construct two functorial factorizations of maps
(see \propref{prop:CofFact} and \propref{prop:TrivCofFact}) that will
be shown in \secref{sec:factor} to be the two factorizations required
by the factorization axiom M5 (see \defref{def:ModCat}).  An important
point for both of these factorizations is that the second map in the
factorization must have the right lifting property with respect to a
set of maps, each of which has a compact domain.  Compact spaces are
``small'' in the sense that, because of
Corollary~\ref{cor:CellGenFiniteSubcmp}, every map from a compact
space into the colimit of a $\lambda$-sequence (see
\defref{def:lambdaseq}) of relative cell complexes must factor through
some intermediate stage of the $\lambda$-sequence, and this is the key
fact that allows us to show that the second map in the factorization
has the required lifting property.  Both of these factorizations are
examples of \emph{the small object argument} (see
\cite{MCATL}*{Proposition~10.5.16}).

\subsection{The first factorization}
\label{sec:FactFirst}

\begin{prop}
  \label{prop:CofFact}
  There is a functorial factorization of every map $f\colon X \to Y$
  as
  \begin{displaymath}
    \xymatrix{
      {X} \ar[r]^{i}
      & {W} \ar[r]^{p}
      & {Y}
    }
  \end{displaymath}
  such that $i$ is a relative cell complex (see \defref{def:CellComp})
  and $p$ has the right lifting property (see \defref{def:LLP}) with
  respect to every element of $I$ (see \defref{def:gc}).
\end{prop}

\begin{proof}
  We will construct the space $W$ as the colimit of a sequence of
  spaces $W_{k}$
  \begin{displaymath}
    \xymatrix{
      {X} \ar@{=}[r]
      & {W_{0}} \ar[r] \ar[d]_{p_{0}}
      & {W_{1}} \ar[r] \ar[dl]_{p_{1}}
      & {W_{2}} \ar[r] \ar[dll]_{p_{2}}
      & {\cdots}\\
      & {Y}
    }
  \end{displaymath}
  where each space $W_{k}$ comes equipped with a map $p_{k}\colon
  W_{k} \to Y$ that makes the diagram commute.  We construct the
  $W_{k}$ inductively, and we begin the induction by letting $W_{0} =
  X$ and letting $p_{0} = f$.

  For the inductive step, we assume that $k \ge 0$ and that we have
  created the diagram through $W_{k}$.  To create $W_{k+1}$, we have
  the solid arrow diagram
  \begin{displaymath}
    \xymatrix@=1em{
      {\coprod S^{n-1}}
      \ar[rr] \ar[dd]
      && {W_{k}} \ar[dd]^{p_{k}} \ar@{..>}[dl]\\
      & {W_{k+1}} \ar@{..>}[dr]^{p_{k+1}}\\
      {\coprod D^{n}}
      \ar[rr] \ar@{..>}[ur]
      && {Y}
    }
  \end{displaymath}
  where the coproducts on the left are indexed by
  \begin{displaymath}
    \coprod_{n \ge 0} \CofSquares \Comma
  \end{displaymath}
  which is the set of commutative squares
  \begin{displaymath}
    \xymatrix{
      {S^{n-1}} \ar[r] \ar[d]
      & {W_{k}} \ar[d]\\
      {D^{n}} \ar[r]
      & {Y \Semicolon}
    }
  \end{displaymath}
  there is one summand for each such square.  We let $W_{k+1}$ be the
  pushout
  \begin{displaymath}
    \xymatrix{
      {\coprod S^{n-1}} \ar[r] \ar[d]
      & {W_{k}} \ar@{..>}[d]\\
      {\coprod D^{n}} \ar@{..>}[r]
      & {W_{k+1} \Period}
    }
  \end{displaymath}
  That defines the space $W_{k+1}$ and the map $p_{k+1}\colon W_{k+1}
  \to Y$.  We let $W = \colim_{k\ge 0} W_{k}$, we let $p\colon W \to
  Y$ be the colimit of the $p_{k}$, and we let $i\colon X \to W$ be
  the composition $X = W_{0} \to \colim W_{k} = W$.

  Since $W$ was constructed by attaching cells to $X$, the map
  $i\colon X \to W$ is a relative cell complex.

  To see that $p\colon W \to Y$ has the right lifting property with
  respect to every element of $I$, suppose that $n \ge 0$ and we have
  the solid arrow diagram
  \begin{equation}
    \label{diag:CofFact}
    \vcenter{
      \xymatrix{
        {S^{n-1}} \ar[r] \ar[d]
        & {W} \ar[d]^{p}\\
        {D^{n}} \ar[r] \ar@{..>}[ur]
        & {Y \Period}
      }
    }
  \end{equation}
  Since $S^{n-1}$ is compact, Corollary~\ref{cor:CellGenFiniteSubcmp} implies
  that there is a positive integer $k$ such that the map $S^{n-1} \to
  W$ factors through $W_{k}$.  Thus, we have the solid arrow diagram
  \begin{displaymath}
    \xymatrix{
      {S^{n-1}} \ar[r] \ar[d]
      & {W_{k}} \ar[r]
      & {W_{k+1}} \ar[r]
      & {W} \ar[d]^{p}\\
      {D^{n}} \ar@{..>}[urr] \ar[rrr]
      &&& {Y}
    }
  \end{displaymath}
  and the map $S^{n-1} \to W_{k}$ is one of the attaching maps in the
  pushout diagram that built $W_{k+1}$ out of $W_{k}$.  Thus, there
  exists a diagonal arrow $D^{n} \to W_{k+1}$ that makes the diagram
  commute, and its composition with $W_{k+1} \to W$ is the diagonal
  arrow required in \eqref{diag:CofFact}.

  To see that the construction is functorial, suppose we have a
  commutative square
  \begin{displaymath}
    \xymatrix{
      {X} \ar[r] \ar[d]_{f}
      & {X'} \ar[d]^{f'}\\
      {Y} \ar[r]
      & {Y' \Semicolon}
    }
  \end{displaymath}
  we will show that the construction of the factorization of $f\colon
  X \to Y$ maps to that of $f'\colon X' \to Y'$, i.e., that there is a
  commutative diagram
  \begin{displaymath}
    \xymatrix@R=3ex{
      &&&&& {Y} \ar[ddd] \\
      {X} \ar@{=}[r]
      & {W_{0}} \ar[r] \ar[d]^{f_{0}} \ar@/^{2ex}/[rrrru]
      & {W_{1}} \ar[r] \ar[d]^{f_{1}} \ar@/^{1ex}/[rrru]
      & {W_{2}} \ar[r] \ar[d]^{f_{2}} \ar[rru]
      & {\cdots} \\
      {X'} \ar@{=}[r]
      & {W_{0}'} \ar[r] \ar@/_{2ex}/[rrrrd]
      & {W_{1}'} \ar[r] \ar@/_{1ex}/[rrrd]
      & {W_{2}'} \ar[r] \ar[rrd]
      & {\cdots} \\
      &&&&& {Y' \Period}
    }
  \end{displaymath}
  We let $f_{0} = f$.  Suppose that we've defined $f_{n}\colon W_{n}
  \to W_{n}'$.  The space $W_{n+1}$ is constructed by attaching an
  $n$-cell to $W_{n}$ for each commutative square
  \begin{displaymath}
    \xymatrix{
      {S^{n-1}} \ar[r]^{\alpha} \ar[d]
      & {W_{n}} \ar[d]\\
      {D^{n}} \ar[r]_{\beta}
      & {Y \Period}
    }
  \end{displaymath}
  We map the cell attached to $W_{n}$ by $\alpha$ to the cell attached
  to $W_{n}'$ by the map $f_{n}\circ \alpha$ indexed by the outer
  commutative rectangle
  \begin{displaymath}
    \xymatrix{
      {S^{n-1}} \ar[r]^{\alpha} \ar[d]
      & {W_{n}} \ar[d] \ar[r]^{f_{n}}
      & {W_{n}'} \ar[d]\\
      {D^{n}} \ar[r]_{\beta}
      & {Y} \ar[r]
      & {Y' \Period}
    }
  \end{displaymath}
  Doing that to each cell attached to $W_{n}$ defines $f_{n+1}\colon
  W_{n+1} \to W_{n+1}'$.
\end{proof}

\subsection{The second factorization}
\label{sec:FactSecond}

\begin{defn}
  \label{def:gtc}
  We will let $J$ denote the set of maps
  \begin{displaymath}
    J = \bigl\{D^{n}\times\{0\} \to D^{n}\times I
        \suchthat n \ge 0\bigr\}
  \end{displaymath}
  and we will call $J$ the \emph{set of generating trivial
    cofibrations}.
\end{defn}

\begin{prop}
  \label{prop:defSfib}
  A map $f\colon X \to Y$ is a Serre fibration if and only if it has
  the right lifting property (see \defref{def:LLP}) with respect to
  every element of $J$ (see \defref{def:gtc}).
\end{prop}

\begin{proof}
  This is just a restatement of the definition of a Serre fibration.
\end{proof}

\begin{defn}
  \label{def:RlJClCmp}
  If $X$ is a subspace of $Y$ such that there is a pushout diagram
  \begin{displaymath}
    \xymatrix{
      {D^{n}\times \{0\}} \ar[r] \ar[d]
      & {X} \ar[d]\\
      {D^{n}\times I} \ar[r]
      & {Y}
    }
  \end{displaymath}
  for some $n \ge 0$, then we will say that $Y$ is obtained from $X$
  by \emph{attaching a $J$-cell} (see \defref{def:gtc}).

  A \emph{relative $J$-cell complex} is an inclusion of a subspace
  $f\colon X \to Y$ such that $Y$ can be constructed from $X$ by a
  (possibly infinite) process of repeatedly attaching $J$-cells.
\end{defn}

\begin{lem}
  \label{lem:Jcell}
  Every element of $J$ is a relative cell complex (see
  \defref{def:CellComp}) with two cells.  If $Y$ is obtained from $X$
  by attaching a $J$-cell, then $X \to Y$ is a relative cell complex
  in which you attach a single $n$-cell and then a single $(n+1)$-cell
  (for some $n \ge 0$).
\end{lem}

\begin{proof}
  There is a homeomorphism between $D^{n}\times I$ and $D^{n+1}$ that
  takes $D^{n}\times\{0\}$ onto one of the two $n$-disks whose union
  is $\partial D^{n+1}$.  Thus, $D^{n+1}$ is homeomorphic to the
  result of first attaching an $n$-cell to $D^{n}\times \{0\}$ and then
  attaching an $(n+1)$-cell to the result, and the pushout of
  \begin{displaymath}
    \xymatrix{
      {D^{n}\times\{0\}} \ar[r] \ar[d]
      & {X}\\
      {D^{n}\times I}
    }
  \end{displaymath}
  is homeomorphic to the result of first attaching an $n$-cell to $X$
  and then attaching an $(n+1)$-cell to the result.
\end{proof}

\begin{rem}
  \label{rem:MultpleJCells}
  We will often construct a relative $J$-cell complex by attaching
  more than one $J$-cell at a time.  That is, given a space $X_{0}$, a
  set $S$, and for each $s \in S$ a map $D^{n_{s}}\times \{0\} \to
  X_{0}$, we may construct a pushout
  \begin{displaymath}
    \xymatrix{
      {\coprod_{s\in S} D^{n_{s}}\times \{0\}} \ar[r] \ar[d]
      & {X_{0}} \ar[d]\\
      {\coprod_{s\in S} D^{n_{s}}\times I} \ar[r]
      & {X_{1}}
    }
  \end{displaymath}
  and then perform a similar construction with $X_{1}$, repeating a
  possibly infinite number of times.
\end{rem}

\begin{prop}
  \label{prop:JClCmpLift}
  Every fibration has the right lifting property (see
  \defref{def:LLP}) with respect to all relative $J$-cell complexes
  and their retracts (see \defref{def:retract}).
\end{prop}

\begin{proof}
  \propref{prop:defSfib} implies that a fibration has the right
  lifting property with respect to every element of $J$, and so
  \lemref{lem:coprodLLP} implies that it has the right lifting
  property with respect to every coproduct of elements of $J$, and so
  \lemref{lem:pushLLP} implies that it has the right lifting property
  with respect to every pushout of a coproduct of elements of $J$, and
  so \lemref{lem:compLLP} implies that it has the right lifting
  property with respect to every composition of pushouts of coproducts
  of elements of $J$, and so \lemref{lem:retractLLP} implies that it
  has the right lifting property with respect to every retract of such
  a composition.
\end{proof}

\begin{prop}
  \label{prop:TrivCofFact}
  There is a functorial factorization of every map $f\colon X \to Y$
  as
  \begin{displaymath}
    \xymatrix{
      {X} \ar[r]^{j}
      & {W} \ar[r]^{q}
      & {Y}
    }
  \end{displaymath}
  such that $j$ is a relative $J$-cell complex (see
  \defref{def:RlJClCmp}) and $q$ has the right lifting property (see
  \defref{def:LLP}) with respect to every element of $J$ (see
  \defref{def:gtc}).
\end{prop}

\begin{proof}
  The construction of the factorization is the same as in
  \propref{prop:CofFact}, except that we use the set $J$ of generating
  trivial cofibrations (see \defref{def:gtc}) in place of the set $I$
  of generating cofibrations (see \defref{def:gc}).

  Since the space $W$ was constructed by attaching $J$-cells to $X$,
  the map $j\colon X \to W$ is a relative $J$-cell complex.

  Since constructing $W_{k+1}$ from $W_{k}$ consists of attaching many
  copies of $D^{n}\times I$ along $D^{n}\times\{0\}$ (for all $n \ge
  0$), \lemref{lem:Jcell} implies that it can be viewed as a 2-step
  process:
  \begin{enumerate}
  \item Attach many $n$-cells (for all $n \ge 0$) to $W_{k}$ to create
    a space we'll call $W'_{k}$.
  \item Attach many $(n+1)$-cells (for all $n \ge 0$) to $W'_{k}$ to
    form $W_{k+1}$.
  \end{enumerate}
  If, for $k \ge 0$, we let $V_{2k} = W_{k}$ and $V_{2k+1} = W'_{k}$,
  then $W$ is the colimit of the sequence
  \begin{displaymath}
    X = V_{0} \to V_{1} \to V_{2} \to \cdots
  \end{displaymath}
  where each space $V_{k+1}$ is built from $V_{k}$ by attaching cells
  to $V_{k}$.  Thus, the map $X \to W$ is a relative cell complex, and
  Corollary~\ref{cor:CellGenFiniteSubcmp} implies that any map
  $D^{n}\times\{0\} \to W$ factors through $W_{k}$ for some $k \ge
  0$.

  The proof that $q\colon W \to Y$ has the right lifting property with
  respect to every element of $J$ and the proof that this construction
  is functorial now proceed exactly as in the proof of
  \propref{prop:CofFact}.
\end{proof}

\begin{prop}
  \label{prop:JclClCmp}
  The relative $J$-cell complex constructed in the proof of
  \propref{prop:TrivCofFact} is a relative cell complex.
\end{prop}

\begin{proof}
  This follows from the proof of \propref{prop:TrivCofFact}.
\end{proof}

\section{The factorization axiom}
\label{sec:factor}

In this section we show that the two factorizations of maps
constructed in \secref{sec:SmObj} are the two factorizations required
by the factorization axiom M5 (see \defref{def:ModCat}).

\subsection{Cofibration and trivial fibration}
\label{sec:CofTrFib}

\begin{prop}
  \label{prop:CofTrFib}
  The functorial factorization constructed in \propref{prop:CofFact}
  of a map $f\colon X \to Y$ as $X \xrightarrow{i} W \xrightarrow{p}
  Y$ where $i$ is a relative cell complex and $p$ has the right
  lifting property with respect to every element of $I$ (see
  \defref{def:gc}) is a factorization into a cofibration followed by a
  map that is both a fibration and a weak equivalence.
\end{prop}

\begin{proof}
  Since the cofibrations are defined to be the relative cell complexes
  and their retracts (see \defref{def:MdCtStr}), the map $i$ is a
  cofibration.

  \lemref{lem:Jcell} implies that every element of $J$ is a relative
  cell complex and \propref{prop:ClCmpLift} implies that the map $p$
  has the right lifting property with respect to all relative cell
  complexes.  Thus, \propref{prop:defSfib} implies that the map $p$ is
  a fibration.

  To see that the map $p$ is a weak equivalence, first note that
  because it has the right lifting property with respect to the map
  $\emptyset \to D^{0}$, it is surjective, and so every path component
  of $Y$ is in the image of a path component of $W$.  Thus, $X$ and
  $Y$ are either both empty or both nonempty.

  To see that every map of homotopy groups (or sets) $\pi_{i}W \to
  \pi_{i}Y$ (for $i \ge 0$) at every basepoint of $W$ is injective,
  note that every element of the kernel gives rise to a commutative
  solid arrow diagram
  \begin{equation}
    \label{diag:CofTrFib}
    \vcenter{
      \xymatrix{
        {S^{i}} \ar[r] \ar[d]
        & {W} \ar[d]^{p}\\
        {D^{i+1}} \ar[r] \ar@{..>}[ur]
        & {Y}
      }
    }
  \end{equation}
  and the existence of the diagonal arrow shows that that element of
  the kernel is the zero element of $\pi_{i}W$.  (If $i = 0$, then the
  fact that this is true for \emph{every} choice of basepoint implies
  that $\pi_{0}W \to \pi_{0}Y$ is a monomorphism.)

  To see that every map of homotopy groups (or sets) $\pi_{i}W \to
  \pi_{i}Y$ (for $i \ge 0$) at every basepoint of $W$ is surjective,
  we first note that since every path component of $Y$ is in the image
  of a path component of $W$, we know that $\pi_{0}W \to \pi_{0}Y$ is
  always surjective.  Now let $i \ge 0$, let $w \in W$, and let $a \in
  \pi_{i+1}\bigl(Y, p(w)\bigr)$.  The element $a$ can be represented
  by a map $D^{i+1} \to Y$ that takes the entire $S^{i}$ that is the
  boundary of $D^{i+1}$ to the point $p(w)$, and so we have a
  commutative solid arrow diagram as in \diagref{diag:CofTrFib} in
  which the upper horizontal map is the constant map to the point $w$
  and the left vertical map is the inclusion.  There must then exist a
  diagonal arrow making the diagram commute, and this is a map
  $D^{i+1} \to W$ that takes the entire boundary $S^{i}$ to the point
  $w$, and thus defines an element of $\pi_{i+1}\bigl(W,w\bigr)$ that
  goes under $p_{*}$ to $a$.
\end{proof}

\subsection{Trivial cofibration and fibration}
\label{sec:TrCofFib}

\begin{prop}
  \label{prop:TrCofFib}
  The functorial factorization constructed in
  \propref{prop:TrivCofFact} of a map $f\colon X \to Y$ as $X
  \xrightarrow{j} W \xrightarrow{q} Y$ where $j$ is a relative
  $J$-cell complex (see \defref{def:RlJClCmp}) and $q$ has the right
  lifting property with respect to every element of $J$ (see
  \defref{def:gtc}) is a factorization into a map that is both a
  cofibration and a weak equivalence followed by a map that is
  fibration.
\end{prop}

\begin{proof}
  The proof of \propref{prop:TrivCofFact} showed that the relative
  $J$-cell complex is a relative cell complex, and so the map $j$ is a
  cofibration.

  Since each inclusion $D^{n}\times\{0\} \to D^{n}\times I$ is the
  inclusion of a strong deformation retract, each map $W_{k} \to
  W_{k+1}$ in the construction of the factorization is also the
  inclusion of a strong deformation retract, and is thus a weak
  equivalence.  Thus, for every $i \ge 0$ the sequence
  \begin{displaymath}
    \pi_{i} W_{0} \to \pi_{i} W_{1} \to \pi_{i} W_{2} \to \cdots
  \end{displaymath}
  is a sequence of isomorphisms.  Since the $S^{i}$ and $D^{i+1}$ are
  all compact, Corollary~\ref{cor:CellGenFiniteSubcmp} implies that every map
  from $S^{i}$ or $D^{i+1}$ to $W$ factors through $W_{k}$ for some $k
  \ge 0$, and so the map $\colim_{k} \pi_{i} W_{k} \to \pi_{i} W$ is
  an isomorphism.  Thus, the map $\pi_{i}X \to \pi_{i}W$ is also an
  isomorphism, and the map $j$ is a weak equivalence.

  \propref{prop:defSfib} implies that the map $q$ is a fibration.
\end{proof}

\section{Homotopy groups and maps of disks}
\label{sec:MpDsk}

The main result of this section is \propref{prop:liftcellwe}, which
will be used in \secref{sec:lift} to prove that the lifting axiom M4
(see \defref{def:ModCat}) holds.  Given a solid arrow diagram as in
\diagref{diag:liftcellwe}, the map $h$ defines an element of
$\pi_{n-1}X$ (at some basepoint), and the existence of the map $g$
implies that the image of that element under $f_{*}$ is the zero
element of $\pi_{n-1}Y$.  Since the map $f$ is a weak equivalence,
this implies that the map $h$ defines the zero element of
$\pi_{n-1}X$, and so $h$ can be extended to a map $G\colon D^{n} \to
X$, but there's no reason for the composition $fG$ to equal $g$.  All
that we know is that $fG$ and $g$ agree on the boundary of $D^{n}$.

We are thus led, in \secref{sec:DiffMp}, to study the situation in
which we have two maps $\alpha,\beta\colon D^{n} \to X$ that agree on
the boundary of $D^{n}$.  We use these maps to define (in
\defref{def:diff}) a \emph{difference map} $d(\alpha,\beta)\colon
S^{n} \to X$, using $\alpha$ on the upper hemisphere and $\beta$ on
the lower hemisphere.  This difference map defines an element of
$\pi_{n}X$ (at some basepoint), and the remaining results of
\secref{sec:DiffMp} show that these homotopy group elements behave as
you would expect.  The results in \secref{sec:LiftDisks} then use the
results of \secref{sec:DiffMp} to replace the map $G$ with one for
which the composition $fG$ equals $g$.

\subsection{Difference maps}
\label{sec:DiffMp}

\begin{defn}
  \label{def:hemisphere}
  For $n \ge 1$ we let $S^{n}_{+}$ be the upper hemisphere of $S^{n}$
  \begin{displaymath}
    S^{n}_{+} = \bigl\{(x_{1}, x_{2}, \ldots, x_{n+1}) \in \R^{n+1}
    \suchthat x_{1}^{2} + x_{2}^{2} + \cdots x_{n+1}^{2} = 1,
    x_{n+1} \ge 0 \bigr\}
  \end{displaymath}
  and we let $S^{n}_{-}$ be the lower hemisphere of $S^{n}$
  \begin{displaymath}
    S^{n}_{-} = \bigl\{(x_{1}, x_{2}, \ldots, x_{n+1}) \in \R^{n+1}
    \suchthat x_{1}^{2} + x_{2}^{2} + \cdots x_{n+1}^{2} = 1,
    x_{n+1} \le 0 \bigr\} \Period
  \end{displaymath}
  We let $p_{+}\colon S^{n}_{+} \to D^{n}$ be the homeomorphism
  \begin{displaymath}
    p_{+}(x_{1}, x_{2}, \ldots, x_{n+1}) =
    (x_{1}, x_{2}, \ldots, x_{n})
  \end{displaymath}
  and we let $p_{-}\colon S^{n}_{-} \to D^{n}$ be the homeomorphism
  \begin{displaymath}
    p_{-}(x_{1}, x_{2}, \ldots, x_{n+1}) =
    (x_{1}, x_{2}, \ldots, x_{n}) \Period
  \end{displaymath}
\end{defn}

\begin{defn}
  \label{def:diff}
  If $X$ is a space and $\alpha,\beta\colon D^{n} \to X$ are maps that
  agree on $\bdry D^{n}$, then we let $d(\alpha,\beta)\colon S^{n} \to
  X$ be the map that is $\alpha \circ p_{+}\colon S^{n}_{+} \to X$ on
  the upper hemisphere of $S^{n}$ and $\beta \circ p_{-}\colon
  S^{n}_{-} \to X$ on the lower hemisphere of $S^{n}$, and we call it
  the \emph{difference map} of $\alpha$ and $\beta$.
\end{defn}

\begin{lem}
  \label{lem:any}
  Let $X$ be a space and let $\alpha\colon D^{n} \to X$ be a map.  If
  $[g] \in \pi_{n}\bigl(X, \alpha(p_{0})\bigr)$ is any element of
  $\pi_{n}\bigl(X, \alpha(p_{0})\bigr)$ (where $p_{0}$ is the
  basepoint of $D^{n}$), then there is a map $\beta\colon D^{n} \to X$
  such that $\beta\rest{\bdry D^{n}} = \alpha\rest{\bdry D^{n}}$ and
  $[d(\alpha,\beta)] = [g]$ in $\pi_{n}\bigl(X, \alpha(p_{0})\bigr)$.
\end{lem}

\begin{proof}
  The basepoint of $D^{n}$ is a strong deformation retract of $D^{n}$,
  and so any two maps $D^{n} \to X$ are homotopic relative to the
  basepoint.  Thus, the restriction of $g$ to $S^{n}_{+}$ is homotopic
  relative to the basepoint to $\alpha \circ p_{+}$.  Since the
  inclusion $S^{n}_{+} \hookrightarrow S^{n}$ has the homotopy
  extension property, there is a homotopy of $g$ to a map $h\colon
  S^{n} \to X$ such that $h\rest{S^{n}_{+}} = \alpha \circ p_{+}$; we
  let $\beta = h \circ (p_{-}^{-1})$, and we have $h =
  d(\alpha,\beta)$.
\end{proof}

\begin{lem}[Additivity of difference maps]
  \label{lem:add}
  If $X$ is a space, $n \ge 1$, and $\alpha,\beta,\gamma\colon D^{n}
  \to X$ are maps that agree on $\bdry D^{n}$, then in
  $\pi_{n}\bigl(X, \alpha(p_{0})\bigr)$ (where $p_{0}$ is the
  basepoint of $D^{n}$) we have
  \begin{displaymath}
    [d(\alpha,\beta)] + [d(\beta,\gamma)] = [d(\alpha,\gamma)]
  \end{displaymath}
  (where, if $n = 1$, addition should be replaced by multiplication).
\end{lem}

\begin{proof}
  Let $T^{n} = S^{n} \cup D^{n}$, where we view $D^{n}$ as the subset
  of $\R^{n+1}$
  \begin{displaymath}
    D^{n} = \bigl\{ (x_{1}, x_{2}, \ldots, x_{n+1}) \suchthat
    x_{1}^{2} + x_{2}^{2} + \cdots x_{n}^{2} \le 1,
    x_{n+1} = 0 \bigr\} \Semicolon
  \end{displaymath}
  $T^{n}$ is then a CW-complex that is the union of the three
  $n$-cells $S^{n}_{+}$, $S^{n}_{-}$, and $D^{n}$, which all share a
  common boundary.  We let $t(\alpha,\beta,\gamma)\colon T^{n} \to X$
  be the map such that
  \begin{align*}
    t(\alpha,\beta,\gamma)\rest{S^{n}_{+}} &= \alpha \circ p_{+}\\
    t(\alpha,\beta,\gamma)\rest{D^{n}} &= \beta\\
    t(\alpha,\beta,\gamma)\rest{S^{n}_{-}} &= \gamma \circ p_{-}
    \Period
  \end{align*}
  We then have that
  \begin{itemize}
  \item the composition $S^{n} \hookrightarrow T^{n}
    \xrightarrow{t(\alpha,\beta,\gamma)} X$ is $d(\alpha,\gamma)$,
  \item the composition $S^{n} \to S^{n}_{+} \cup D^{n} \subset T^{n}
    \xrightarrow{t(\alpha,\beta,\gamma)} X$ (where that first map is
    the identity on $S^{n}_{+}$ and is $p_{-}$ on $S^{n}_{-}$) is
    $d(\alpha,\beta)$, and
  \item the composition $S^{n} \to D^{n} \cup S^{n}_{-} \subset T^{n}
    \xrightarrow{t(\alpha,\beta,\gamma)} X$ (where that first map is
    $p_{+}$ on $S^{n}_{+}$ and is the identity on $S^{n}_{-}$) is
    $d(\beta,\gamma)$.
  \end{itemize}

  The basepoint of $D^{n}$ is a strong deformation retract of $D^{n}$,
  and so the map $\beta\colon D^{n} \to X$ is homotopic relative to
  the basepoint to the constant map to $\alpha(p_{0})$ (where $p_{0}$
  is the common basepoint of $S^{n}$ and $D^{n}$).  Since the
  inclusion $D^{n} \hookrightarrow T^{n}$ has the homotopy extension
  property, there is a homotopy of $t(\alpha,\beta,\gamma)$ relative
  to the basepoint to a map $\hat t(\alpha,\beta,\gamma)$ that takes
  all of $D^{n}$ to the basepoint $\alpha(p_{0})$, and
  $d(\alpha,\gamma)$ is homotopic to the composition
  \begin{displaymath}
    S^{n} \hookrightarrow T^{n}
    \xrightarrow{\hat t(\alpha,\beta,\gamma)} X \Period
  \end{displaymath}

  If $T^{n} \to S^{n} \vee S^{n}$ is the map that collapses $D^{n}$ to
  a point, then $\hat t(\alpha,\beta,\gamma)$ factors as $T^{n} \to
  S^{n} \vee S^{n} \xrightarrow{\alpha_{\beta} \vee \beta_{\gamma}}
  X$, where $\alpha_{\beta}\colon S^{n} \to X$ is homotopic to
  $d(\alpha,\beta)$ and $\beta_{\gamma}\colon S^{n} \to X$ is
  homotopic to $d(\beta,\gamma)$.  Thus, $d(\alpha,\gamma)$ is
  homotopic to the composition
  \begin{displaymath}
    S^{n} \hookrightarrow T^{n} \longrightarrow S^{n} \vee S^{n}
    \xrightarrow{\alpha_{\beta} \vee \beta_{\gamma}} X \Period
  \end{displaymath}
  and so we have $[d(\alpha,\gamma)] = [d(\alpha,\beta)] +
  [d(\beta,\gamma)]$ if $n>1$ and $[d(\alpha,\gamma)] =
  [d(\alpha,\beta)] \cdot [d(\beta,\gamma)]$ if $n=1$.
\end{proof}

\begin{lem}
  \label{lem:DskHmtp}
  If $X$ is a space, $\alpha,\beta\colon D^{n} \to X$ are maps that
  agree on $\bdry D^{n}$, and $[d(\alpha,\beta)]$ is the identity
  element of $\pi_{n}X$, then $\alpha$ and $\beta$ are homotopic
  relative to $\bdry D^{n}$.
\end{lem}

\begin{proof}
  Since $[d(\alpha,\beta)]$ is the identity element of $\pi_{n}X$,
  there is a map $h\colon D^{n+1} \to X$ whose restriction to $\bdry
  D^{n+1}$ is $d(\alpha,\beta)$.  View $D^{n}\times I$ as the cone on
  $\bdry(D^{n}\times I) = (D^{n}\times \{0\}) \cup (S^{n-1}\times I)
  \cup (D^{n}\times\{1\})$ with vertex at the center of $D^{n}\times
  I$.  Let $p\colon D^{n}\times I \to D^{n+1}$ be the map that
  \begin{itemize}
  \item on $D^{n}\times \{0\}$ is the composition $D^{n}\times\{0\}
    \xrightarrow{\pr} D^{n} \xrightarrow{(p_{+})^{-1}} S^{n}_{+}
    \hookrightarrow D^{n+1}$,
  \item on $D^{n}\times \{1\}$ is the composition $D^{n}\times\{0\}
    \xrightarrow{\pr} D^{n} \xrightarrow{(p_{-})^{-1}} S^{n}_{-}
    \hookrightarrow D^{n+1}$,
  \item on $S^{n-1}\times I$ is the composition $S^{n-1}\times I
    \xrightarrow{\pr} S^{n-1} \hookrightarrow S^{n}_{+}\cap S^{n}_{-}
    \subset D^{n+1}$,
  \item takes the center point of $D^{n}\times I$ to the center point
    of $D^{n+1}$, and
  \item is linear on each straight line connecting the center point of
    $D^{n}\times I$ to its boundary.
  \end{itemize}
    The composition $D^{n}\times I \xrightarrow{p} D^{n+1}
    \xrightarrow{h} X$ is then a homotopy from $\alpha$ to $\beta$
    relative to $\bdry D^{n}$.
\end{proof}

\subsection{Lifting maps of disks}
\label{sec:LiftDisks}

In this section, we use the results of \secref{sec:DiffMp} to prove
\propref{prop:liftcellwe}.  Given a solid arrow diagram as in
\diagref{diag:liftcellwe}, \propref{prop:lift} shows that, even if the
weak equivalence $f\colon X \to Y$ isn't a fibration, there exists a
diagonal arrow that makes the upper triangle commute exactly and makes
the lower triangle commute up to a homotopy relative to the boundary
of $D^{n}$.  \propref{prop:liftcell} then shows that if the weak
equivalence $f$ is also fibration, then there exists a diagonal arrow
that makes both triangles commute exactly.  Both \propref{prop:lift}
and \propref{prop:liftcell} actually only apply to one value of $n$ at
a time, and \propref{prop:liftcellwe} is the statement that if
\emph{all} homotopy groups (and sets) go by an isomorphism under
$f_{*}$, then the diagonal arrow exists for all values of $n$.

\begin{prop}
  \label{prop:lift}
  Let $f\colon X \to Y$ be a map, let $n \ge 1$, and suppose that we
  have the solid arrow diagram
  \begin{displaymath}
    \xymatrix{
      {\bdry D^{n}} \ar[r]^-{h} \ar[d]_{i}
      & {X} \ar[d]^{f}\\
      {D^{n}} \ar[r]_-{g} \ar@{..>}[ur]
      & {Y \Period}
    }
  \end{displaymath}
  If $F$ is the homotopy fiber of $f$ over some point in the image of
  $g$ and if $\pi_{n-1}F = 0$, then there exists a map $G\colon D^{n}
  \to X$ such that $Gi = h$ and $fG \homotopic g$ relative to
  $\bdry D^{n}$.
\end{prop}

\begin{proof}
  The map $h$ defines an element $[h]$ of $\pi_{n-1}X$ (at some
  basepoint) such that $f_{*}\bigl([h]\bigr) = 0$ in $\pi_{n-1}Y$.
  Since $\pi_{n-1}F = 0$, the long exact homotopy sequence of a
  fibration implies that $[h] = 0$ in $\pi_{n-1}X$, and so there is a
  map $j\colon D^{n} \to X$ such that $j\circ i = h$.

  The maps $fj\colon D^{n} \to Y$ and $g\colon D^{n} \to Y$ agree on
  $\partial D^{n}$, and so there is a difference map $d(fj,g)\colon
  S^{n} \to X$ (see \defref{def:diff}) that defines an element
  $\alpha$ of $\pi_{n}Y$.  Since $\pi_{n-1}F = 0$, the long exact
  homotopy sequence implies that there is an element $\beta$ of
  $\pi_{n}X$ such that $f_{*}(\beta) = -\alpha$ (if $n>1$) or
  $f_{*}(\beta) = \alpha^{-1}$ (if $n=1$), and \lemref{lem:any}
  implies that we can choose a map $G\colon D^{n} \to X$ that agrees
  with $j\colon D^{n} \to X$ on $\partial D^{n}$ such that $[d(G,j)] =
  \beta$ in $\pi_{n}X$.  Thus, $Gi = h$, and since $[d(fG, fj)] =
  [f\circ d(G,j)] = f_{*}[d(G,j)] = f_{*}(\beta) = -\alpha$ (if $n>1$)
  or $\alpha^{-1}$ (if $n=1$), \lemref{lem:add} implies that
  $[d(fG,g)] = [d(fG,fj)] + [d(fj,g)] = -\alpha + \alpha = 0$ (with a
  similar statement if $n=1$), and so \lemref{lem:DskHmtp} implies
  that $fG$ is homotopic to $g$ relative to $\bdry D^{n}$.
\end{proof}

\begin{lem}
  \label{lem:fibr}
  A map $f\colon X \to Y$ is a Serre fibration if and only if for
  every $n \ge 0$ and every solid arrow diagram
  \begin{equation}
    \label{diag:fibrOth}
    \vcenter{
      \xymatrix{
        {(D^{n}\times\{0\}) \cup (\partial D^{n}\times I)}
        \ar[r] \ar[d]
        & {X} \ar[d]^{f}\\
        {D^{n}\times I} \ar[r] \ar@{..>}[ur]
        & {Y}
      }
    }
  \end{equation}
  there exists a diagonal arrow making both triangles commute.
\end{lem}

\begin{proof}
  For every $n \ge 0$ there is a homeomorphism of pairs
  \begin{displaymath}
    \bigl(D^{n}\times I, (D^{n}\times\{0\}) \cup
    (\partial D^{n}\times I)\bigr)
    \longrightarrow
    \bigl(D^{n}\times I, D^{n}\times\{0\}\bigr)
  \end{displaymath}
  under which diagrams of the form \eqref{diag:fibrOth} correspond to
  diagrams of the form
  \begin{displaymath}
    \xymatrix{
      {D^{n}\times \{0\}} \ar[r] \ar[d]
      & {X} \ar[d]^{f}\\
      {D^{n}\times I} \ar[r] \ar@{..>}[ur]
      & {Y}
    }
  \end{displaymath}
  and \propref{prop:defSfib} implies that there always exists a
  diagonal arrow making the diagram commute if and only if $f$ is a
  fibration.
\end{proof}

\begin{prop}
  \label{prop:liftcell}
  Let $f\colon X \to Y$ be a fibration, let $n \ge 1$, and suppose
  that we have the solid arrow diagram
  \begin{displaymath}
    \xymatrix{
      {\partial D^{n}} \ar[r]^-{h} \ar[d]_{i}
      & {X} \ar[d]^{f}\\
      {D^{n}} \ar[r]_-{g} \ar@{..>}[ur]
      & {Y}
    }
  \end{displaymath}
  If $F$ is the fiber of $f$ over some point in the image of $g$ and
  if $\pi_{n-1}F = 0$, then there exists a diagonal arrow making both
  triangles commute.
\end{prop}

\begin{proof}
  \propref{prop:lift} implies that there is a map $G\colon D^{n} \to
  X$ such that $Gi = h$ and $fG \homotopic g$ relative to $\bdry
  D^{n}$.  Let $H\colon D^{n}\times I \to Y$ be a homotopy from $fG$
  to $g$ relative to $\partial D^{n}$.  We have a lift to $X$ of the
  restriction of $H$ to $(D^{n}\times\{0\}) \cup (\partial D^{n}\times
  I)$ defined as $G \circ \pr_{D^{n}}$ on $D^{n}\times \{0\}$ and $h
  \circ \pr_{\partial D^{n}}$ on $\partial D^{n}\times I$, and
  \lemref{lem:fibr} implies that we can lift the homotopy $H$ to a
  homotopy $H'\colon D^{n}\times I \to X$ such that the restriction of
  $H'$ to $D^{n}\times\{0\}$ is $G \circ \pr_{D^{n}}$ and the
  restriction of $H'$ to $\partial D^{n}\times I$ is $h \circ
  \pr_{\partial D^{n}}$.  Let $G'\colon D^{n} \to X$ be defined by
  $G'(d) = H'(d,1)$, and we have $G' \circ i = h$ and $f\circ G' = g$.
\end{proof}

\begin{prop}
  \label{prop:liftcellwe}
  If $f\colon X \to Y$ is both a fibration and a weak equivalence,
  then for every $n \ge 0$ and every solid arrow diagram
  \begin{equation}
    \label{diag:liftcellwe}
    \vcenter{
      \xymatrix{
        {\partial D^{n}} \ar[r]^-{h} \ar[d]_{i}
        & {X} \ar[d]^{f}\\
        {D^{n}} \ar[r]_-{g} \ar@{..>}[ur]
        & {Y}
      }
    }
  \end{equation}
  there exists a diagonal arrow making the diagram commute.
\end{prop}

\begin{proof}
  Since $f$ is a weak equivalence, the set of path components of $X$
  maps onto that of $Y$, and so an application of
  \propref{prop:defSfib} with $n = 0$ implies that $X$ maps onto $Y$.
  That implies the case $n = 0$, and the cases $n \ge 1$ follow from
  \propref{prop:liftcell}.
\end{proof}

\section{The lifting axiom}
\label{sec:lift}

In this section we use the results of \secref{sec:MpDsk} (in
particular, \propref{prop:liftcellwe}) to show that the lifting axiom
M4 (see \defref{def:ModCat}) is satisfied (see \thmref{thm:onelift}
and \thmref{thm:otherlift}).

\subsection{Cofibrations and trivial fibrations}
\label{sec:CofTFib}

\begin{thm}
  \label{thm:onelift}
  If $f\colon X \to Y$ is both a fibration and a weak equivalence,
  then it has the right lifting property with respect to all
  cofibrations.
\end{thm}

\begin{proof}
  \propref{prop:liftcellwe} implies that $f$ has the right lifting
  property with respect to every element of $I$ (see \defref{def:gc}).
  Since a cofibration is defined as a retract of a relative cell
  complex, the result follows from \propref{prop:ClCmpLift}.
\end{proof}

\subsection{Trivial cofibrations and fibrations}
\label{sec:FibTCof}

\begin{thm}
  \label{thm:otherlift}
  If $i\colon A \to B$ is both a cofibration and a weak equivalence,
  then it has the left lifting property with respect to all
  fibrations.
\end{thm}

\begin{proof}
  \propref{prop:TrivCofFact} implies that we can factor $i\colon A \to
  B$ as $A\xrightarrow{s} W \xrightarrow{t} B$ where $s$ is a relative
  $J$-cell complex and $t$ is a fibration (see
  \propref{prop:defSfib}).  \propref{prop:TrCofFib} implies that $s$
  is both a cofibration and a weak equivalence and $t$ is a fibration.
  Since $s$ and $i$ are weak equivalences, the ``two out of three''
  property of weak equivalences implies that $t$ is also a weak
  equivalence, and so \thmref{thm:onelift} implies that $i$ has the
  left lifting property with respect to $t$.  The retract argument
  (see \propref{prop:RetArg}) now implies that $i$ is a retract of
  $s$, and so \propref{prop:JClCmpLift} implies that $i$ has the left
  lifting property with respect to all fibrations.
\end{proof}

\begin{cor}
  \label{cor:RelCWlift}
  If $i\colon A \to B$ is both a relative CW-complex and a weak
  equivalence and $p\colon X \to Y$ is a Serre fibration, then for
  every solid arrow diagram
  \begin{displaymath}
    \xymatrix{
      {A} \ar[r] \ar[d]_{i}
      & {X} \ar[d]^{f}\\
      {B} \ar[r] \ar@{..>}[ur]
      & {Y}
    }
  \end{displaymath}
  there is a dotted arrow making the diagram commute.
\end{cor}

\begin{proof}
  Since every relative CW-complex is a relative cell complex, the map
  $i$ is both a cofibration (see \defref{def:MdCtStr}) and a weak
  equivalence, and since our fibrations are Serre fibrations (see
  \defref{def:MdCtStr}), the result now follows from
  \thmref{thm:otherlift}.
\end{proof}

\section{The proof of \thmref{thm:main}}
\label{sec:proofmain}

We must show that the five axioms of \defref{def:ModCat} are satisfied
by the weak equivalences, cofibrations, and fibrations of
\defref{def:MdCtStr}.

The limit axiom M1 is satisfied because we've assumed that our
category of spaces is complete and cocomplete.

For the two out of three axiom M2, we first note that if any two of
$f$, $g$, and $gf$ induce an isomorphism of the set of path
components, then so does the third.  If our maps are $f\colon X \to Y$
and $g\colon Y \to Z$, then if either $f$ and $g$ are weak
equivalences or $g$ and $gf$ are weak equivalences, then the two out
of three property applied to homomorphisms of homotopy groups implies
that the third map also induces an isomorphism of homotopy groups at
an arbitrary choice of basepoint.  If $f$ and $gf$ are assumed to be
weak equivalences, then we know that every choice of basepoint in $Y$
is in the same path component as a point in the image of $f$, and so
(using the change of basepoint isomorphism) it is sufficient to show
that the homotopy groups at a basepoint in the image of $f$ are mapped
isomorphically, and that follows from the two out of three property of
group homomorphisms.

For the retract axiom M3,
\begin{itemize}
\item a retract of a weak equivalence is a weak equivalence because a
  retract of a group isomorphism (or of an isomorphism of the set of
  path components) is an isomorphism,
\item a retract of a cofibration is a cofibration because a retract of
  a retract of a relative cell complex is a retract of a relative cell
  complex, and
\item a retract of a fibration is a fibration because of
  \propref{prop:defSfib} and \lemref{lem:retractLLP}.
\end{itemize}

For the lifting axiom M4, when $i$ is a cofibration and $p$ is both a
fibration and a weak equivalence, this is \thmref{thm:onelift}, and
when $i$ is both a cofibration and a weak equivalence and $p$ is a
fibration, this is \thmref{thm:otherlift}.

For the factorization axiom M5, the factorization into a cofibration
followed by a map that is both a fibration and a weak equivalence is
\propref{prop:CofTrFib}, and the factorization into a map that is both
a cofibration and a weak equivalence followed by a fibration is
\propref{prop:TrCofFib}.



\end{document}